\title{Strong approximation for the supermarket model}
\documentclass[12pt]{article}

\usepackage{amssymb}
\usepackage{latexsym}
\usepackage{amsmath}
\usepackage{amssymb}
\usepackage{amsthm}

\newtheorem{theorem}{Theorem}[section]

\newtheorem{proposition}[theorem]{Proposition}

\newcommand{\halmos}{\rule{1ex}{1.4ex}}
\newcommand{\proofbox}{\hspace*{\fill}\mbox{$\halmos$}}

\newcommand{\R}{\mathbb{R}}
\newcommand{\N}{\mathbb{N}}
\newcommand{\nats}{\mathbb{N}}

\renewcommand{\a}{\alpha}
\renewcommand{\d}{\delta}
\newcommand{\s}{\sigma}
\renewcommand{\O}{\Omega}
\newcommand{\half}{\frac12}

\newcommand{\E}{\mathbb{E}}
\newcommand{\pr}{\mathbb{P}}

\newcommand{\p}{\mathbb{P}}

\author{Malwina J Luczak
\thanks{Department of Mathematics, London School of Economics, Houghton Street, London WC2A 2AE, UK.  e-mail: malwina@planck.lse.ac.uk}
\& James Norris
\thanks{Statistical Laboratory, Centre for Mathematical Sciences, Wilberforce Road, Cambridge CB3 0WB, UK. e-mail: J.R.Norris@statslab.cam.ac.uk}}

\begin{document}

\maketitle

\begin{abstract}
We prove three strong approximation theorems for the ``supermarket'' or
``join the shortest queue'' model -- a law of large numbers, a jump process
approximation and a central limit theorem. The estimates are carried through
rather explicitly. This allows us to estimate each of the infinitely many
components of the process in its own scale and to exhibit a cut-off in the
set of active components which grows slowly with the number of servers.

\textit{Key words and phrases.} Supermarket model, join the shortest queue, 
law of large numbers, diffusion approximation, exponential martingale 
inequalities.

\textit{MSC (2000).} Primary: 60K25 . Secondary: 60F15. 
\end{abstract}

\section{Introduction}
\label{sec.intro}

The supermarket model is a system of $N$ single-server queues. Customers
arrive as a Poisson process of rate $N \lambda$, where $\lambda \in (0,1)$. 
Each customer examines $d$ queues, chosen randomly from all $N$ queues, where
$d \ge 2$, and joins the shortest of these $d$ queues, choosing randomly
if the shortest queue is not unique. The service times of all customers are
independent rate 1 exponential random variables. We will be concerned with 
the behaviour of this model when $\lambda$ and $d$ are fixed, over a finite
time interval $[0,t_0]$, as $N \rightarrow \infty$. We shall consider the case
when the system starts in some well-behaved state with low server loads (in a 
sense to be made precise below).

This model has attracted attention because it turns out that the choice
offered to customers, even if $d=2$, dramatically reduces queue lengths; 
see~\cite{G99,M96,VDK96}; and in particular, the length of the longest queue,
see~\cite{LM03,LM04}. Given that our analysis relies on
$N$ being very much larger than $d$, the model does not describe well the 
behaviour of a real supermarket. Rather it serves as an example where a simple
dynamic routing rule leads to a greatly improved performance, which is
of interest in the context of communications networks.

Our results provide strong approximations for the supermarket model, 
and include
a law of large numbers and a diffusion approximation. In arriving at these
results we have developed techniques to establish weak convergence of 
a sequence of Markov processes $X^N$ in infinitely many dimensions, 
where the jumps of $X^N$ are of order $N^{-1}$ and occur at a rate of order $N$.
The classical results for fluid limits are set in a finite-dimensional
context. We make essential use of the fact that the number of ``active'' 
components in $X^N$ grows
only very slowly with $N$. We have used direct and quantitative methods based
on exponential martingales and strong approximation of Poisson processes by
Brownian motion. These methods seem well-suited to deal with such ``almost
finite-dimensional'' Markov processes. Earlier results for this model include
laws of large numbers in~\cite{G99,G00,M96,VDK96}, quantitative concentration
of measure estimates in~\cite{LM04} and a central limit 
theorem~\cite{G04}. See also~\cite{L03}
for a preliminary version of the law of large numbers presented in the present
paper (for a more general range of initial conditions).

The limiting behaviour of the supermarket model as $N \rightarrow \infty$ may
conveniently be described in terms of the vector $X_t=(X^{k}_t: k \in 
\nats)$, where $X^{k}_t$ denotes the proportion of all $N$ queues having at 
least $k$ customers at time $t$. The process $X=(X_t)_{t \ge 0}$ is a Markov 
chain. We will suppose throughout that $X_0=x_0$ with $x_0$ non-random and we
will suppress the dependence of $X$ on $N$ to lighten the notation. Now $X$
has the form of a density dependent Markov chain such as considered by
Ethier and Kurtz in~\cite{EK86}, Chapter 11. Thus one might expect to be able 
to find a deterministic process $(x_t)_{t \ge 0}$ and a Gaussian process 
$(\gamma_t)_{t \ge 0}$ such that
\[X_t=x_t+O(N^{-1/2}), \quad X_t =x_t +N^{-1/2} \gamma_t +O(\log N/N).\]
However the number of non-zero components in $X$ grows with $N$
so the standard theory
does not apply. 
%Moreover the second, Gaussian, approximation relies on a non-degeneracy in 
%the limiting diffusion which may fail even in finite dimensions.

We will see that, for small initial data,
the component $X^{k}$ has a scale $a_k=\lambda^{1+d+\ldots
+ d^{k-1}}$, which of course decays very rapidly with $k$. Thus the number of
queues having at least $k$ customers is of order $Na_k$. We can find $m$ of 
order $\log \log N$ such that $Na_m$ is of order 1. 
Thus we can exhibit a cut-off in 
the number of active components which grows only slowly with $N$. 
Below the cut-off, for $k\le m-1$, we prove convergence with explicit
control of error probabilities for each of the $\log \log N$ active components. 

We thereby obtain results of the form
\begin{eqnarray*}
& & X^{k}_t =x^k_t + a_k O(\sqrt{\log \log \log N/Na_k} ),\\
& & X^{k}_t = x^k_t + N^{-1/2}\gamma^k_t + 
a_k O(\log {(Na_k)}/Na_k).
\end{eqnarray*}
Note that each component is estimated in the correct scale, with an error
depending on the number of queues active at that level. The $\log \log \log N$
in the first equation is a (small) price we pay for working with infinitely 
many components.
These asymptotics will be established with a degree of uniformity in $x_0$, 
which thus allows a dependence of $x_0$ on $N$.
The Gaussian approximation relies, as in the finite-dimensional case, on a 
sophisticated coupling of the compensated Poisson process with Brownian motion
due to Koml\'os, Major and Tusn\'ady~\cite{KMT75}. 

We will give a third result, for $k\le m-1$, of the form
\[X^{k}_t = x^k_t  + N^{-1/2}\tilde{\gamma}^{k}_t +
a_k O( { (\log \log \log N/Na_k)}^{3/4} ).\]
Here $(\tilde{\gamma}_t)_{t \ge 0}$ is a jump process with drift which depends
on $N$ but is of a simpler type than $X$ in that it is a 
linear function of additive Poisson noise.
The characteristics of $\tilde{\gamma}$ are derived in a simple and 
canonical way from those of $X$. 
Moreover $\tilde{\gamma}$ and $X$ share a common filtration. 
The error term is larger than for the Gaussian approximation. On the other hand
the derivation is significantly simpler.

We obtain also the behaviour of the queue sizes at and above the cut-off.
We see a residual randomness in 
$(X^{m}_t)_{t \ge 0}$ even for large values of $N$. This may be approximated
in terms of an $M/M/\infty$ queue with arrival
rate $N\lambda(x^{m-1}_t)^d$ and service rate 1. Over a 
given finite time interval, there are no queues with lengths greater than $m$.

Thus we will show for the supermarket model that its infinite-dimensional 
character does not prevent the derivation of precise asymptotics. We expect the
general approach taken here to adapt well to a number of further examples
of similar character.

\section{Statement of results}
\label{sec.results}

Let $S_0$ denote the set of non-increasing sequences $x=(x^k: k \in \nats )$ 
in $[0,1]$. For $x\in S_0$, set $x^0=1$ and
define $\lambda_{+}^k(x) =\lambda ({(x^{k-1})}^d-{(x^k)}^d)$,
$\lambda_{-}^k (x)=x^k-x^{k+1}$ and  
$b^k(x)=\lambda_{+}^k(x) -\lambda_{-}^k(x)$. 
It is shown in \cite{VDK96} that, given $x_0\in S_0$, there is a unique
solution $(x_t)_{t\ge0}$ to $\dot x_t=b(x_t)$ in $S_0$.
Moreover, for any other solution $(y_t)_{t\ge0}$ in $S_0$,
$x^k_0\le y^k_0$ for all $k$ implies $x^k_t\le y^k_t$ for all $k$ 
and all $t\ge0$.

Recall that $a_k=\lambda^{1+d+ \ldots + d^{k-1}}$ for $k \in \nats$.
Then $a=(a_k: k \in \nats)$ is the unique solution in $S_0$ to $b(x)=0$ 
such that $\lim_{k \to \infty} x^k =0$; 
%which vanishes at $\infty$; 
to see this note that $x^{k+1}-\lambda {(x^k)}^d$ is
independent of $k$ for any solution $x$.
Define $\|x\|=\sup_k|x^k|/a_k$, set $E=\{x\in \mathbb R^{\nats}:\|x\|<\infty\}$
and set $S=S_0\cap E$. 
Extend $\lambda_\pm$ and $b$ to $\R^\N$ by setting 
$\lambda_{+}^k(x) =\lambda ({(y^{k-1})}^d-{(y^k)}^d)^+$,
$\lambda_{-}^k (x)=(y^k-y^{k+1})^+$, where $y^k=(x^k)^+$.
%\m{Extension to $E$ needed for next sentence.}
It is easy to check that $b$ maps $E$ to itself and is locally
Lipschitz for the given norm. Thus, if $x_0\in S$, then $(x_t)_{t\ge0}$
does not leave $S$ immediately.
Moreover, if $\|x_0\|\le1$, then, by comparison with the stationary 
solution $a$, $\|x_t\|\le1$ for all $t\ge0$.
For $\rho\ge1$ and $t_0>0$ set 
\[ S(\rho,t_0)=\{x_0\in S:\|x_t\|\le\rho\mbox{ for all }t\le t_0\}.\]

The state-space $I$ of the Markov chain $X=(X_t)_{t \ge 0}$ is the set of
non-increasing sequences in $N^{-1} {\{0,1,\ldots , N\}}$ with finitely many
non-zero terms. Thus $I\subseteq S$. The L\'evy kernel for $X$ is given by
\[ K(x,dy)=\sum_{k=1}^{\infty} \left [ N \lambda_{+}^k (x) \delta_{e_k/N}(dy)
+N\lambda_{-}^k(x) \delta_{-e_k/N}(dy) \right ],\]
where $e_k$ denotes the $k$-th standard basis vector.
Given $m\in\nats$,
let $(\hat{X}^{m}_t)_{t\ge t_0}$ be a process starting from $x_0^m$ and such
that $N\hat{X}^m$ is an $M/M/\infty$ queue,
with arrival rate $N\lambda (x^{m-1}_t)^d$ and service
rate $1$.  
We can now state our law of large numbers.

\begin{theorem}
\label{thm.fluid-limit}
Set $m=m(N)=\inf\{k\in\nats: Na_k\le (\log N)^4\}$.
There is a coupling of $\hat X^m$ and $X^m$ such that, 
for all $\rho\ge1,t_0>0$ and all sequences 
$R(N)$ with $R(N)/\sqrt{\log\log\log N} \to\infty$, we have
\begin{align*}
\sup_{x_0\in I\cap S(\rho,t_0)}
\pr_{x_0} (&\sqrt{N}|X^{N,k}_t-x^k_t| > R(N) \sqrt{a_k} 
\mbox{ for some }k \le m-1\\
&\mbox{  or  } X^{N,m}_t\not = \hat{X}^{N,m}_t \mbox{  or  }
X^{N,m+1}_t \not = 0 \mbox{  for some  } t \le t_0)\rightarrow 0.
\end{align*} 
In particular
\[\sup_{k \le m-1} \sup_{t \le t_0} |X^{k,N}_t-x^k_t|/a_k \rightarrow 0\]
in probability, uniformly in $x_0\in I\cap S(\rho,t_0)$.
\end{theorem}
We know, see~\cite{VDK96}, that if $\rho > 0$ and $\|x_0\|\le \rho$ 
then $x^k_t \rightarrow a_k$ as $t \rightarrow \infty$.
Thus, for $k\le m-1$,  the proportional
error in approximating $X^{k}_t$ by the deterministic process $x^k_t$ is 
small for large values of $N$.

The central limit theorem shows generically that the power $\sqrt N$ in
Theorem~\ref{thm.fluid-limit} cannot be improved while the approximating
process $(x_t)_{t \ge 0}$ remains deterministic. Our next result is a
refined approximation which allows an improvement to $N^{3/4}$.
Let $\tilde{\mu}$ be a Poisson random measure on 
$\mathbb R^{\mathbb N}\times(0,t_0]$ with intensity
\[\tilde\nu(dy,dt)=K(x_t,dy)dt.\]
We show in Section~\ref{App} that the linear equations 
\begin{eqnarray}
\label{eq.jump.approx}
\tilde{\gamma}_t^k=\sqrt{N}\int_{\mathbb R^{\mathbb N}\times(0,t]}
y^k(\tilde{\mu} -\tilde{\nu})(dy,ds) + \int_0^t \nabla b^k(x_s)
\tilde{\gamma}_s ds
\end{eqnarray}
have a unique cadlag solution $(\tilde{\gamma}_t^k:k\in\nats,t\le t_0)$.
Set $\tilde X_t=x_t+N^{-1/2}\tilde{\gamma}_t$.

\begin{theorem}
\label{thm.jump.approx}
Define $m(N)$ as in Theorem \ref{thm.fluid-limit}. 
There is a coupling of $\tilde X$ and $X$, in a common filtration, such that, 
for all $\rho\ge1,t_0>0$ and all
sequences $\tilde R(N)$ with $\tilde R(N)/(\log\log\log N)^{3/4}\to\infty$,
we have
\begin{equation*}
\sup_{x_0\in I\cap S(\rho,t_0)}
\pr_{x_0} (N^{3/4}|X^k_t-\tilde X^k_t| > \tilde R(N)a_k^{1/4}
\mbox{ for some }k \le m-1, t \le t_0)\rightarrow 0.
\end{equation*}
\end{theorem}

The final result is a diffusion approximation. 
Let $B^k_+, B^k_-, k \in \nats$, be independent standard Brownian motions.
Set $\sigma^k_{\pm}(x)=\sqrt{\lambda^k_{\pm}(x)}$. 
We show in Section~\ref{App} that the linear equations 
\begin{equation}
\label{eq.diffusion}
\gamma^k_t=\int_0^t \sigma^k_+(x_s) dB^k_+(s)
-\int_0^t \sigma^k_-(x_s) dB^k_-(s)+ \int_0^t \nabla b^k(x_s)
\gamma_s ds,\quad t\le t_0
\end{equation}
have a unique solution $(\gamma_t^k:k\in\nats,t\le t_0)$ with
$$
\sup_{k\in\nats}\E(\sup_{t\le t_0}|\gamma_t^k|^2)<\infty.
$$ 
Set $\bar X_t=x_t+N^{-1/2}\gamma_t$.

\begin{theorem}
\label{thm.diffusion}
Define $m(N)$ as in Theorem \ref{thm.fluid-limit}.           
There is a coupling of $\bar X$ and $X$ such that,
for all $\rho\ge1,t_0>0$, there is a constant $\bar R$, independent of $N$,
such that
\begin{equation*}
\sup_{x_0\in I\cap S(\rho,t_0)}
\pr_{x_0} (N|X^k_t-\bar X^k_t| > \bar R\log(Na_k)
\mbox{ for some }k \le m-1, t \le t_0)\rightarrow 0.
\end{equation*}
\end{theorem}

We remark that there are alternative versions of all three theorems
in which $x_0$ is replaced by $x_0^{(m)}=(x_0^1,\dots,x_0^{m-1},0,\dots)$ 
and $\lambda_\pm(x)$
is replaced by $\lambda_\pm(x)^{(m)}$, so that the approximating 
deterministic dynamics are ${(m-1)}$-dimensional.
The proofs are a minor modification of the proofs given below.
These alternative versions would have merit in any computational 
implementation of the approximations since $m$
is of order only $\log\log N$.

\section{Law of large numbers}
\label{sec.lln}

In the first half of this section, we fix $N$ and $A,R\ge1$ and set
$m=\inf\{k\in\nats:Na_k\le A\}$. We will
obtain, subject to certain constraints, a global estimate on the 
probability appearing in Theorem~\ref{thm.fluid-limit}. In the second half
we will show that this estimate becomes small as $N\to\infty$ when
$A=(\log N)^4$ and when $R$ is chosen as in Theorem~\ref{thm.fluid-limit}.

Fix $\rho\ge1$ and $t_0>0$ and assume that $\rho A^d\le N^{d-1}$.
Set $\alpha = (\log \log N-\log \log (1/\lambda ))/\log d$, so that
$N \lambda^{d^\alpha}=1$ and $Na_k=\lambda^{1+d+ \ldots + d^{k-1}-d^{\alpha}}$.
If $k \le \alpha$ then $N a_k \ge 1$, whereas if $k \ge \alpha +1$ then $N a_k
\le 1$. Hence, at least for sufficiently large $N$, we will have $m \in (\alpha
-1,\alpha +2)$.

Consider the case
$x_0\in S(\rho, t_0)$. Then $Nx^{m+1}_0\le N\rho a_{m+1}<N\rho a_m^d
\le\rho A^d/N^{d-1}\le1$, so $x_0^{m+1}=0$.
Set $T_1=\inf\{t\ge0:X^{m+1}_t\not=0\}$.
Note that, while $X^{m+1}_t=0$, $X^{m+1}_t$ increases at rate 
$N\lambda (X^m_{t-})^d$, whereas $X^m_t+X^{m+1}_t$ increases at rate 
$N\lambda (X^{m-1}_{t-})^d$ and decreases at rate $NX^m_{t-}$.
We can therefore find an $M/M/\infty$ queue $(Q_t)_{t\ge0}$, 
starting from $Nx_0^m$, without arrivals 
and with service rate $1$, and a Poisson random measure 
$\mu(dt,dx,du)$ on $(0,\infty)^3$, independent of $Q$
and of intensity $e^{-u}dtdxdu$, such that, for $t\le T_1$,
\[NX^{m+1}_t=\mu(\{(s,x,u):s\le t<s+u,x\le N\lambda (X^m_{s-})^d\})\]
and
\[N(X^m_t+X^{m+1}_t)=Q_t
+\mu(\{(s,x,u):s\le t<s+u,x\le N\lambda (X^{m-1}_{s-})^d\}).\]
Define $(\hat X^m_t)_{t\ge0}$ by
\[N\hat X^m_t=Q_t
+\mu(\{(s,x,u):s\le t<s+u,x\le N\lambda (x^{m-1}_s)^d\})\]
and set $T_2=\inf\{t\ge0:X^m_t+X^{m+1}_t\not=\hat X^m_t\}$.
Then $(N\hat X^m_t)_{t\ge0}$ is an $M/M/\infty$ queue starting from
$Nx_0^m$, with arrival rate
$N\lambda(x^{m-1}_t)^d$ and service rate $1$.
Fix $r>1$ and set $T_3=\inf\{t\ge0:\hat X^m_t>ra_m\}$.
Fix $R\ge1$, set
\[T_4=\inf\{t\ge0:\sqrt N|X^k_t-x_t^k|>R\sqrt a_k\mbox{ for some }k\le m-1\},\]
and set $T=T_1\wedge T_2\wedge T_3\wedge T_4\wedge t_0$. Finally, set
\[p=p(N,\lambda, d, x_0,A,R,r)=\p(T<t_0).\]

\begin{proposition}\label{PP}
Assume that $x_0\in S(\rho,t_0)$, that $A,R,\rho\ge1$ 
with $\rho A^d\le N^{d-1}$, that $r>\rho$ and that 
\[2rAt_0e^{Lt_0}/N^{(1/2)(1-1/d)}\le R\le(t_0\wedge1)\sqrt{A},\] 
where $L=2(d\s^{d-1}+1)$ and $\s=\rho+1$. 
Then $p\le p_1+p_2+p_3+p_4$, where
\begin{align*}
p_1= & A^dr^dt_0/N^{d-1},\\
p_2= & A^{1-1/(2d)}d\s^{d-1}Rt_0/N^{(1/2)(1-1/d)},\\
%p_2= & \sqrt Ndr^{d-1}Ra_mt_0/\sqrt{a_{m-1}},\\
p_3= & \rho^dt_0/(r-\rho),\\
p_4= & 2m\exp (-R^2/(20\s^dt_0e^{2Lt_0})).
\end{align*}
\end{proposition}
\begin{proof}
It will suffice to show that $\p(T=T_i)\le p_i$ for $i=1,2,3,4$.
Recall that for a Poisson random variable $Y$ of
parameter $\nu>0$ and for $a>0$ we have $\p(Y\ge a)\le\nu/a$.
For $t<T$ we have $N\lambda(X^m_t)^d\le N\lambda(ra_m)^d\le A^dr^d/N^{d-1}$,
so $X^{m+1}_T$ is dominated by a Poisson random variable $Y_1$ of
parameter $p_1$, and so 
\[\p(T=T_1)=\p(X^{m+1}_T=1)\le\p(Y_1\ge1)\le p_1.\]
Since $x_0\in S(\rho,t_0)$, we have
$x_t^k\le \rho a_k$ for all $k\in\nats$ and all $t\le t_0$.
For $k\le m-1$ we have $R\sqrt{a_k/N}\le a_kR/\sqrt{Na_k}
\le a_kR/\sqrt A\le a_k$.
Hence, for $t<T$ and $k\le m-1$,
\begin{equation}\label{Xa}
X^k_t\le x_t^k+R\sqrt{a_k/N}\le \s a_k.
\end{equation}
Then, for $t<T$,
\begin{align*}
N\lambda|(X^{m-1}_t)^d & -(x^{m-1}_t)^d|
 \le N\lambda d\s^{d-1}a_{m-1}^{d-1}|X^{m-1}_t-x^{m-1}_t|\\
& \le N\lambda d\s^{d-1}a_{m-1}^{d-1}R\sqrt{a_{m-1}/N}
\le A^{1-1/(2d)}d\s^{d-1}R/N^{(1/2)(1-1/d)}.
\end{align*}
Set 
\[\Delta=\mu(\{(t,x,u):t\le T, N\lambda(x^{m-1}_t\wedge X^{m-1}_{t-})^d
<x\le N\lambda(x^{m-1}_t\vee X^{m-1}_{t-})^d\}).\]
Then $\Delta$ is dominated by a Poisson random variable $Y_2$ of
parameter $p_2$. Hence
\[\p(T=T_2)=\p(\Delta=1)\le\p(Y_2\ge 1) \le p_2.\]
Note that $\lambda(x^{m-1}_t)^d\le\lambda\rho^da_{m-1}^d=\rho^da_m$ 
for all $t\le t_0$.
Thus $N\hat X^m_T\le Nx^m_0+Y_3$ for a Poisson random variable $Y_3$
of parameter $N\rho^da_mt_0$ and so
\[\p(T=T_3)=\p(\hat X^m_T\ge ra_m)\le\p(Y_3\ge N(r-\rho)a_m)\le p_3.\]
It remains to estimate $\p(T=T_4)$. For this we write
\begin{equation}\label{X}
X_t^k=x_0^k+M_t^k+\int_0^tb^k(X_s)ds
\end{equation}
so that 
\[X_t^k-x_t^k=M_t^k+\int_0^t(b^k(X_s)-b^k(x_s))ds.\]
Then we use a combination of exponential martingale inequalities and Gronwall's 
lemma to obtain the desired estimate. 
First we investigate how small $M$ will need to be to obtain the
required bound on $|X_t^k-x_t^k|$ for $k\le m-1$. Note that
$$
|\lambda_+^k(x)-\lambda_+^k(y)|
\le  \lambda d(x^{k-1} \lor y^{k-1})^{d-1}|x^{k-1} -y^{k-1}| 
+  \lambda d(x^k \lor y^k)^{d-1}|x^k-y^k|
$$
so, provided that $x^k \lor y^k \le\s a_k$ and $x^{k-1} \lor y^{k-1}
\le \s a_{k-1}$,
\begin{equation}\label{B+}
|\lambda_+^k(x)-\lambda_+^k(y)|
\le   d\s^{d-1}\{(a_k/a_{k-1})|x^{k-1} -y^{k-1}| + |x^k-y^k|\}.
\end{equation}
Also
\begin{equation}\label{B-}
|\lambda_-^k(x)-\lambda_-^k(y)|
\le|x^k-y^k| + |x^{k+1}-y^{k+1}|.
\end{equation}
Hence, provided that $x^k \lor y^k \le\s a_k$ and $x^{k-1} \lor y^{k-1}
\le \s a_{k-1}$, we have
$$
|b^k(x)-b^k(y)|/\sqrt{a_k} 
 \le  L \sup_{j =k-1,k,k+1}|x^j-y^j|/\sqrt{a_j}.
$$
We note that the definitions of $T_2$ and $T_3$ force $X^m_t\le ra_m$ 
for all $t<T$.
Set
\[f(t)=\sup_{k\le m-1}\sup_{s\le t} |X_s^k-x_s^k|/\sqrt{a_k}.\]
Then, for $t<T$ and $k\le m-2$,
\[|b^k(X_t)-b^k(x_t)|/\sqrt{a_k}\le Lf(t)\]
and
\[|b^{m-1}(X_t)-b^{m-1}(x_t)|/\sqrt{a_{m-1}}\le Lf(t)+ra_m/\sqrt{a_{m-1}}.\]
Hence, for $t\le T$,
\[f(t)\le(M^*_t+ra_mt/\sqrt{a_{m-1}})+L\int_0^tf(s)ds\]
where
\[M^*_t=\sup_{k\le m-1}\sup_{s\le t} |M^k_s|/\sqrt{a_k}.\]
Set $\a_k=\half e^{-Lt_0}R\sqrt{a_k/N}$ and consider for $k\le m-1$ 
the stopping times $T^k=T^k_-\wedge T^k_+$, where
\[T^k_\pm=\inf\{t\ge0:\pm M_t>\a_k\}.\]
Suppose that $T<T^1\wedge\dots\wedge T^{m-1}$.
Then $M^*_T\le\half e^{-Lt_0}R/\sqrt{N}$.
On the other hand
\[ra_mt_0/\sqrt{a_{m-1}}=\lambda^{1/(2d)}ra_m^{1-1/(2d)}t_0
\le rA^{1-1/(2d)}t_0/N^{1-1/(2d)}\le(1/2)e^{-Lt_0}R/N^{1/2}.\]
So by Gronwall's lemma
\[ f(T)\le e^{Lt_0}(M^*_T+ra_mt_0/\sqrt{a_{m-1}})\le R/\sqrt{N}.\]
Hence
\[\p(T=T_4)\le\sum_{k=1}^{m-1}\p(T^k\le T)\]
and it remains to estimate $\p(T^k_\pm\le T)$ for $k\le m-1$.
For $k\in\nats,x\in S$ and $\theta\in\R$, set 
\[ \phi^k(x,\theta)=\lambda_+^k(x)h(\theta )+\lambda_-^k(x) h(-\theta),\]
where $h(\theta )=e^{\theta}-1-\theta$. For $t<T$ and $k\le m-1$ we have
$X^k_t\le \s a_k$ so
\begin{equation}\label{PH}
\phi^k(X_t,\theta)\le\lambda \s^da_{k-1}^d h(\theta )+\s a_kh(-\theta)
\le \s^da_kg(\theta )
\end{equation}
where $g(\theta )=e^\theta -2+e^{-\theta}$.
Consider for $\theta\ge0$  the exponential martingale
\begin{align*}
Z^k_t&=\exp\{N\theta(X^k_t-X^k_0)-\int_0^t\int_{\mathbb R^{\mathbb N}}
(e^{N\theta y^k}-1)K(X_s,dy)ds\}\\
     &=\exp {\{N \theta M^k_t-N \int_0^t {\phi^k (X_s, \theta)ds}\}}
\end{align*}
and note that on the event $T^k_+ \le T$ we have
\[Z^k_{T^k_+}\ge\exp {\{N \theta \alpha_k - N\s^da_kg(\theta)t_0\}}.\]
By optional stopping $\E (Z^{k,N}_{T^k_1})\le1$, so by Chebyshev's inequality
\[ \pr (T^k_+ \le T)\le \exp (-N\theta\a_k+N\s^da_kg(\theta)t_0).\]
We choose $\theta=\a_k/(2\s^da_k\gamma t_0)$, where $\gamma=g(1)\le5/4$.
It is straightforward to check that $\theta\le1$ so 
$g(\theta)\le\gamma\theta^2$. Hence
\[\pr (T^k_+ \le T)\le \exp (-N\a_k^2/(4\s^da_k\gamma t_0))
=\exp (-R^2/(16\s^d\gamma t_0e^{2Lt_0})).\]
The same bound applies to $\pr (T^k_- \le T)$. So we have shown that
$\p(T_4=T)\le p_4$ as required.
\end{proof}

\begin{proof}[Proof of Theorem \ref{thm.fluid-limit}]
We will determine conditions on sequences $R(N)$ and $r(N)$ so that,
for $A(N)=(\log N)^4$, 
as $N\to\infty$, all the constraints of Proposition \ref{PP} are satisfied
and, with an obvious notation, $p_i(N)\to0$ for $i=1,2,3,4$.
For $p_4(N)\to0$ it suffices that 
$\log\log N\exp(-R^2/(20\s^dt_0e^{2Lt_0}))\to0$
and hence that $R/\sqrt{\log\log\log N}\to\infty$.
For $p_3(n)\to0$ it suffices that $r\to\infty$.
For $p_2(N)\to0$ it suffices that $A^{2-1/d}R^2/N^{1-1/d}\to0$
and for $p_1(N)\to0$ it suffices that $Ar/N^{1-1/d}\to0$.
If we can also arrange that $R/\sqrt{A}\to0$ and $rA/(RN^{(1/2)(1-1/d)})\to0$
then all the constraints of Proposition \ref{PP} will be satisfied eventually.
A possible choice is to take 
$r(N)=N^{(1/2)(1-1/d)}/(\log N)^4$ and any sequence $R(N)$ with
$R(N)/\sqrt{\log\log\log N}\to\infty$
and $R(N)/(\log N)^2\to0$.
This proves the first part of the theorem.
For the final assertion it suffices to note that for $k\le m-1$,
$R(N)/\sqrt{Na_k}\le R(N)/\sqrt{A(N)}\to0$.
\end{proof}

We remark that the choice $R(N)=\log N$ leads to a bound of the form
$p(N)\le C N^{-(1/2)(1-1/d)}$ up to logarithmic corrections. This is the best
rate of decay of probabilities we have found. We remark also that a marginally
shorter proof can be had by replacing the exponential martingale inequality
by Doob's $L^2$-inequality, at the small cost of requiring that
$R(N)/\sqrt{\log\log N}\to\infty$.

%%%%%%%%%%%%%%%%%%%%%%%%%%%%%%%%%%%%%%%%%%%%%%%%%%%%%%%%%%%%%%%%%%%%%%%%%%%%%%%%%%%%%%%%
%%%%%%%%%%%%%%%%%%%%%%%%%%%%%%%%%%%%%%%%%%%%%%%%%%%%%%%%%%%%%%%%%%%%%%%%%%%%%%%%%%%%%%%%

\section{A refinement of the fluid limit}
\label{sec.jump.approx}

This section leads to a proof of Theorem \ref{thm.jump.approx}.
The deterministic limit (for components $k\le m-1$) just discussed will be 
refined by approximating the martingale $M$ in \eqref{X} by another
martingale whose characteristics are determined by the limit path,
and at the same time linearizing around the limit path. The accuracy of the
approximation is thereby improved from $N^{-1/2}$ to $N^{-3/4}$ at the cost
of moving to an approximating process which is not deterministic but has
a simple random structure, being a linear function of a Poisson random
measure.

Define a measure $\tilde\nu$ on $\mathbb R^{\mathbb N}\times (0,t_0]$ by 
\[\tilde\nu(dy,dt)=K(x_t,dy)dt.\]
We will take $\tilde\mu$ to be a Poisson random measure with intensity
$\tilde\nu$ coupled, in a way to be specified, with the process $X$.
Define $\tilde M=(\tilde M^k_t:k\in\nats, t\le t_0)$ by
\[\tilde M_t^k=\int_{\mathbb R^{\mathbb N}\times(0,t]}y^k(\tilde\mu-\tilde\nu)(dy,ds)\]
and define $\tilde\gamma=(\tilde\gamma^k_t:k\in\nats, t\le t_0)$ by
\[\tilde\gamma_t=\sqrt{N}\tilde M_t+\int_0^t\nabla b(x_s)\tilde\gamma_sds.\]
We show in Section~\ref{App} 
that we can write $\tilde\gamma$ as an explicit linear function
of $\tilde\mu-\tilde\nu$ 
\begin{equation}\label{GF}
\tilde\gamma_t=\sqrt{N}\int_{\mathbb R^{\mathbb N}\times(0,t]}
\Phi_{t,s}y(\tilde\mu-\tilde\nu)(dy,ds)
\end{equation}
where $(\Phi_{t,s}:s\le t\le t_0)$ is the $\nats\times\nats$ matrix-valued process
given by 
\[\frac\partial{\partial t}\Phi_{t,s}=\nabla b(x_t)\Phi_{t,s},
\quad \Phi_{s,s}=I.\]
Thus $\tilde\gamma$ has a simpler stochastic structure than $X$.
In particular,
we can write the characteristic function of any finite-dimensional
distribution of $\tilde\gamma$ in terms of $(x_t)_{t\le t_0}$ 
and $(\Phi_{t,s}:s\le t\le t_0)$.

Recall that
\[X_t=x_0+M_t+\int_0^tb(X_s)ds.\]
On the other hand, if we set $\tilde X_t=x_t+N^{-1/2}\tilde\gamma_t$, then
\[\tilde X_t=x_0+\tilde M_t+\int_0^tb(x_s)ds+\int_0^t\nabla b(x_s)
(\tilde X_s-x_s)ds.\]
Set $\tilde Y=X-\tilde X,\tilde D=M-\tilde M$ and
\begin{equation}\label{A}
A_t=\int_0^t(b(X_s)-b(x_s)-\nabla b(x_s)(X_s-x_s))ds.
\end{equation}
Then
\begin{equation}\label{Y}
\tilde Y_t=\tilde D_t+A_t+\int_0^t\nabla b(x_s)\tilde Y_sds.
\end{equation}
We will obtain a good approximation if we can couple $\tilde M$ with $M$ 
to make $\tilde D$ small. Define kernels
$K_0,K_+, K_-$ on $(0,t_0]\times E\times \mathbb R^{\mathbb N}$ by
\begin{eqnarray*}
&&K_0(t,x,dy)=K(x,dy) \land K(x_t,dy),\\
&&K_{\pm}(t,x,dy)=(K(x,dy)-K(x_t,dy))^{\pm},
\end{eqnarray*}
and let $K_*(t,x,w,dx',dw')$ be the image of the measure
\[K_0(t,x,dy_0) \otimes K_+(t,x,dy_+ ) \otimes K_-(t,x,dy_- )\]
by the map $(x',w')=(y_0+y_+,y_0+y_-)$. Let $(X_t,W_t)_{t \ge 0}$ 
be a Markov chain, starting from $(x_0,0)$, with time-dependent
L\'evy kernel $K_*$. Set
\[\tilde{\mu}_t = \sum_{\Delta W_t \not = 0}\delta_{(t,\Delta W_t)}.\]
Then $(X_t)_{t \ge 0}$ is a Markov chain with L\'evy kernel $K$
and $\tilde\mu$ is a Poisson random measure with intensity
$\tilde\nu$. 
We have coupled $X$ and $W$ so that, as far as possible, they have
the same jumps. 
Set $T_5=\inf\{t\ge0:|\tilde\gamma_t^m|>r\sqrt{a_m}\}$, where $r$ is as in 
the previous section. 
Fix $\tilde R>0$ and set 
\[T_6=\inf\{t\ge0:N^{3/4}|X_t^k-\tilde X_t^k|>\tilde Ra_k^{1/4}
\mbox{ for some }k\le m-1\}.\]
Finally, set $\tilde T=T\wedge T_5\wedge T_6$ and $\tilde p= \p(\tilde T<t_0)$.

\begin{proposition}\label{QQ}
Assume that the conditions of Proposition~\ref{PP} hold. 
Set ${H=\frac12d(d-1)\s^{d-2}}$
and assume in addition that $rA\le N^{(1/2)(1-1/d)}$ and
\[4HR^2t_0e^{Lt_0}/N^{1/4}+4rAt_0e^{Lt_0}/N^{(1/4)(1-1/d)}
\le\tilde R
\le 4RLt_0e^{Lt_0}A^{1/4}.\]
Then $\tilde p\le p_1+p_2+p_3+p_4+p_5+p_6$, where $p_1,p_2,p_3,p_4$ are
defined in Proposition~\ref{PP} and
\[p_5=8(\rho^d+1)t_0e^{2Lt_0^2}/r^2, 
\quad p_6=2m\exp\{-\tilde R^2/(20RLt_0e^{2Lt_0})\}.\]
\end{proposition}
\begin{proof}
Given Proposition~\ref{PP}, 
it will suffice to show that $\p(\tilde T=T_i)\le p_i$ for $i=5,6$.
By Proposition~\ref{GT}, 
$$
\E(\sup_{t\le t_0}|\tilde\gamma_t^m|^2)\le8(\rho^d+1)t_0e^{2Lt_0^2}a_m.
$$
So $\p(\tilde T=T_5)\le p_5$ by Chebyshev's inequality.

We now follow an argument similar to the proof that $\p(T=T_4)\le p_4$
in Proposition \ref{PP}.
Set
\begin{align*}
\tilde f(t)= & \sup_{k\le m-1}\sup_{s\le t}|\tilde Y_s^k|/a_k^{1/4},\\
\tilde A^*_t= & \sup_{k\le m-1}\sup_{s\le t}|A^k_s|/a_k^{1/4}
+\int_0^t|\tilde Y_s^m|ds/a_{m-1}^{1/4},\\
\tilde D^*_t= & \sup_{k\le m-1}\sup_{s\le t}|\tilde D_s^k|/a_k^{1/4}.
\end{align*}
We note that
\[\nabla b^k(x)y=\lambda d(x^{k-1})^{d-1}y^{k-1}-\lambda d(x^k)^{d-1}y^k
-y^k+y^{k+1}\]
so, provided that $x^k\le \s a_k$ and $x^{k-1}\le\s a_{k-1}$,
\[|\nabla b^k(x)y|/a_k^{1/4}\le L\sup_{j=k-1,k,k+1}|y^j|/a_j^{1/4}.\]
For $t\le t_0$ we have $x_t^k\le\rho a_k\le\s a_k$ for all $k$, 
so for $k\le m-1$
\[|\nabla b^k(x_t)\tilde Y_t|/a_k^{1/4}\le L\tilde f(t)
+\delta_{k,m-1}|\tilde Y^m_t|/a_{m-1}^{1/4}.\]
Then, from equation (\ref{Y}) we get 
\[\tilde f(t)\le\tilde D^*_t+\tilde A^*_t+L\int_0^t\tilde f(s)ds,\]
so, by Gronwall's lemma, $\tilde f(t)\le e^{Lt}(\tilde A^*_t+\tilde D^*_t)$ 
for all $t\le t_0$.
Note that, for $k\le m-1$,
\begin{align*}
b^k(y) & -b^k(x)-\nabla b^k(x)(y-x)\\
& =\lambda((y^{k-1})^d-(x^{k-1})^d-d(x^{k-1})^{d-1}(y^{k-1}-x^{k-1}))\\
& \quad-\lambda((y^k)^d-(x^k)^d-d(x^k)^{d-1}(y^k-x^k))
\end{align*}
so, provided that $x^k,y^k\le \s a_k$ and $x^{k-1},y^{k-1}\le\s a_{k-1}$,
\begin{align*}
|b^k(y) & -b^k(x)-\nabla b^k(x)(y-x)|\\
& \le H\lambda(a_{k-1}^{d-2}|y^{k-1}-x^{k-1}|^2
  +a_k^{d-2}|y^k-x^k|^2 ).
\end{align*}
For $t<\tilde T$ and $k\le m-1$ we have $|X_t^k-x^k_t|\le R\sqrt{a_k/N}$; 
moreover, as we showed at (\ref{Xa}), this implies that $X_t^k\le\s a_k$.
Hence
\begin{align*}
|b^k(X_t) & -b^k(x_t)-\nabla b^k(x_t)(X_t-x_t)|/a_k^{1/4}\\
& \le H\lambda R^2(a_{k-1}^{d-1} +a_k^{d-1})/(Na_k^{1/4})\le 2CR^2/N.
\end{align*}
Also, $|\tilde Y_t^m|\le|X^m_t-x_t^m|+N^{-1/2}|\tilde\gamma^m_t|$ so, 
for $t<\tilde T$, $|\tilde Y_t^m|/a_{m-1}^{1/4}
\le(ra_m+N^{-1/2}r\sqrt{a_m})/a_{m-1}^{1/4}
\le 2r(A/N)^{1-1/(4d)}$.
It follows that 
$$
\tilde A^*_{\tilde T}\le 2HR^2t_0/N+2rt_0(A/N)^{1-1/(4d)}
\le \half e^{-Lt_0}\tilde R/N^{3/4}.
$$
Set
$\tilde\alpha_k=\frac12 e^{-Lt_0}\tilde Ra_k^{1/4}/N^{3/4}$ and consider
the stopping times $\tilde T^k=\tilde T^k_+\wedge\tilde T^k_-$, where
\[\tilde T^k_\pm=\inf\{t\ge0:\pm \tilde D_t^k>\tilde\alpha_k\}.\]
Suppose that ${\tilde T}<\tilde T^1\wedge\dots\wedge\tilde T^{m-1}$. Then
$\tilde D^*_{\tilde T}\le\half e^{-Lt_0}\tilde R/N^{3/4}$
so $f({\tilde T})\le \tilde R/N^{3/4}$ and ${\tilde T}<T_6$.
Hence
\[\p(\tilde T=T_6)\le\sum_{k=1}^{m-1}\p(\tilde T^k\le\tilde T)\]
and it remains to estimate $\p(\tilde T^k_\pm\le\tilde T)$ for $k\le m-1$.

For $k\le m-1$ set
\begin{align*}
\psi^k (t,x,\theta ) = 
& (\lambda^k_+(x)-\lambda^k_+(x_t))^+ h(\theta)
+ (\lambda^k_-(x)-\lambda^k_-(x_t))^+ h(-\theta )\\
& +(\lambda^k_+(x)-\lambda^k_+(x_t))^- h(-\theta)
+ (\lambda^k_-(x)-\lambda^k_-(x_t))^- h(\theta ).
\end{align*}
Fix $\theta\ge0$ and consider for $k\le m-1$ the exponential martingale
\begin{align*}
\tilde{Z}^k_t
& = \exp\{N\theta(X_t^k-X_0^k-W_t^k)
    -\int_0^t\int_{\mathbb R^{\mathbb N}}(e^{N\theta y^k}-1)K^+(s,X_s,dy)ds\\
& \quad\quad\quad\quad\quad\quad 
    -\int_0^t\int_{\mathbb R^{\mathbb N}}(e^{-N\theta y^k}-1)K^-(s,X_s,dy)ds\}\\
& =\exp\{N\theta \tilde D_t^k-N\int_0^t\psi^k (s,X_s,\theta )ds\}.
\end{align*}
For $k\le m-2$ and $t\le t_0$, for $x^k\le\s a_k$ and $x^{k-1}\le\s a_{k-1}$, 
we can estimate as at~(\ref{B+}),~(\ref{B-}),~(\ref{PH}) to obtain
\begin{align*}
\psi^k (t,x,\theta )/\sqrt{a_k}
&\le g(\theta)(|\lambda^k_+(x)-\lambda^k_+(x_t)|
  +|\lambda^k_-(x)-\lambda^k_-(x_t)|)/\sqrt{a_k}\\
& \le Lg(\theta)\sup_{j=k-1,k,k+1}|x^j-x_t^j|/\sqrt{a_j}
\end{align*}
so, for $t<\tilde T$ and $k\le m-2$,
\begin{equation}\label{PS}
\psi^k (t,X_t,\theta )\le Lg(\theta) R\sqrt{a_k/N}.
\end{equation}
Similarly, since we assume $R\ge1,rA\le N^{(1/2)(1-1/d)}$ and $Na_m\le A$, 
we have, for $t<\tilde T$, $X^m_t\le ra_m\le R\sqrt{a_{m-1}/N}$ and we can
show that (\ref{PS}) remains true for $k=m-1$. By optional stopping
we have $\mathbb E(\tilde{Z}^k_{\tilde T^k_+})\le1$ for all $k\le m-1$. But on 
the event $\tilde T^k_+\le\tilde T$ we have
$$\tilde{Z}^k_{\tilde T^k_+}\ge\exp\{N\theta\tilde\a^k-Lg(\theta)
Rt_0\sqrt{Na_k}\}.$$
We choose $\theta=2\sqrt{N}\tilde\a^k/(5LRt_0\sqrt{a_k})$, checking that
$\theta\le1$, so that $g(\theta)\le5\theta^2/4$, 
and use Chebyshev's inequality to deduce
$$\p(\tilde T^k_+\le\tilde T)\le\exp\{-N^{3/2}\tilde\a_k^2/(5LRt_0\sqrt{a_k})\}
=\exp\{-\tilde R^2/(20LRt_0e^{2Lt_0})\}.
$$
The same bound applies to $\pr (\tilde T^k_- \le\tilde T)$. So we have shown that
$\p(\tilde T=T_6)\le p_6$ as required.
\end{proof}

\begin{proof}[Proof of Theorem \ref{thm.jump.approx}]
Choose $r(N)$ as in the proof of Theorem~\ref{thm.fluid-limit}, so
that $r(N)\to\infty$ and so $p_5(N)\to0$. Assume $\tilde{R}(N)/(\log N)^2 \to
0$, and set 
$s(N)=\tilde R(N)/(\log\log\log N)^{3/4}$. Then $s(N)\to\infty$.
Set $R(N)=s(N)(\log\log\log N)^{1/2}$ and $r(N)=N^{(1/4)(1-1/d)}/(\log N)^4$.
It is straightforward to check that all the constraints in Propositions
\ref{PP} and \ref{QQ} are satisfied eventually. Moreover, as in the proof
of Theorem \ref{thm.fluid-limit}, we have $p_i(N)\to0$ for $i=1,2,3,4$.
Finally $\tilde R(N)^2/R\log\log\log N\to\infty$ so also $p_6(N)\to0$
which proves the theorem.
\end{proof}

%%%%%%%%%%%%%%%%%%%%%%%%%%%%%%%%%%%%%%%%%%%%%%%%%%%%%%%%%%%%%%%%%%%%%%%%%%%%%%%%%%%%%%
%%%%%%%%%%%%%%%%%%%%%%%%%%%%%%%%%%%%%%%%%%%%%%%%%%%%%%%%%%%%%%%%%%%%%%%%%%%%%%%%%%%%%%%

\section{Diffusion approximation}
\label{sec.diffusion}

In this section we prove Theorem~\ref{thm.diffusion}. 
The method follows the lines set out in \cite{EK86}, Chapter 11.
As we have already seen, our process $X$ has around $\log\log N$
active components, which have a wide range of scales.
This will require special consideration in the implementation of the
general method.
We also have to deal with the fact that the variance of the diffusion
approximation has degeneracies. The diffusion coefficient, obtained
as the square root of the variance, then fails to be Lipschitz and some
special care is needed to arrive at the desired convergence.

Let $(X^k_t:k\in\nats,t\ge0)$ be the supermarket process starting from $x_0$
and recall equation (\ref{X})
\begin{equation*}
X_t^k=x_0^k+M_t^k+\int_0^tb^k(X_s)ds.
\end{equation*}
Recall also that we set $\bar X_t=x_t+N^{-1/2}\gamma_t$, where 
$(\gamma_t^k:k\in\nats,t\le t_0)$ is defined by the linear equations
(\ref{eq.diffusion})
\begin{equation*}
\gamma^k_t=\sqrt N\bar M^k_t+ \int_0^t \nabla b^k(x_s) \gamma_s ds.
\end{equation*}
and
\begin{equation*}
\sqrt N\bar M^k_t=\int_0^t\sigma^k_+(x_s)dB^k_+(s)-\int_0^t\sigma^k_-(x_s) dB^k_-(s).
\end{equation*}
Set $Y=X-\bar X$ and $D=M-\bar M$. Then
$$
Y_t=D_t+A_t+\int_0^t \nabla b^k(x_s) Y_s ds
$$
where $A_t$ is defined at (\ref{A}).
We will obtain a good approximation if we can couple $\bar M$ with $M$ 
to make $D$ small. 

The coupling relies on the following approximation result of \cite{KMT75}:
there exists a constant $c\in(0,\infty)$ and a probability space on which are
defined a compensated Poisson process $Z$ of rate $1$ and a standard
Brownian motion $W$ such that, for all $t\ge0$ and $x\in\R$,
\begin{equation}\label{KM}
\p(\sup_{s\le t}|Z(s)-W(s)|\ge c\log t+x)\le ce^{-x/c}.
\end{equation}

Given independent compensated Poisson processes $Z_+^k,Z_-^k,k\in\nats$ of 
rate $1$ we can construct $X$ by the equations (\ref{X}) and
$$
M_t^k=N^{-1}\{Z_+^k(N\int_0^t\lambda^k_+(X_s)ds)-
Z_-^k(N\int_0^t\lambda^k_-(X_s)ds)\}.
$$
On the other hand, by a theorem of Knight, see for example \cite{RY91}, 
there exist independent Brownian
motions $W_+^k,W_-^k,k\in\nats$, such that, for all $k\in\nats$ and $t\le t_0$,
$$
W_\pm^k(N\int_0^t\lambda_\pm^k(\bar X_s)ds)
=\sqrt N\int_0^t\s_\pm^k(\bar X_s)dB_\pm^k(s).
$$
The law of $(B_+^k,B_-^k:k\in\nats)$ given $(W_+^k,W_-^k:k\in\nats)$ is given by
a measurable kernel. So we may assume that these processes are defined on the 
same probability space as $(Z_+^k,Z_-^k:k\in\nats)$ and that
$(Z_+^k,W_+^k),(Z_-^k,W_-^k),k\in\nats,$ are independent copies of $(Z,W)$.

Set $T_7=\inf\{t\ge0:|\gamma_t^m|>r\sqrt{a_m}\}$.
Fix $\bar R>0$ and set
\[T_8=\inf\{t\ge0:N|X_t^k-\bar X_t^k|>\bar R\log(Na_k)
\mbox{ for some }k\le m-1\}.\]
Finally, set $\bar T=T\wedge T_7\wedge T_8$ and $\bar p= \p(\bar T<t_0)$.

\begin{proposition}\label{RR}
Assume that the conditions of Proposition~\ref{PP} hold, and assume moreover
that $A \ge e^2$ and $R\le\sqrt A/2$.
There is a constant $C<\infty$, depending only on $d,\lambda, \rho$ and $t_0$
such that, if 
$$C+C(R^2+rA)/\log N \le\bar R\le A/(2\log A)$$ 
then $\bar p\le p_1+p_2+p_3+p_4+p_7+p_8$, where $p_1,p_2,p_3,p_4$ are
defined in Proposition~\ref{PP} and
%\[p_7=4\rho^d t_0e^{Lt_0}/r^2, \quad p_8=\dots.\]
\[p_7=C/r^2, \quad p_8=Cm(A^{-1}+(\bar R\log A)^{-2}).\]
\end{proposition}
\begin{proof}
Given Proposition~\ref{PP},
it will suffice to show that $\p(\bar T=T_i)\le p_i$ for $i=7,8$.
By Proposition~\ref{GT},
$$
\E(\sup_{t\le t_0}|\gamma_t^m|^2)\le8(\rho^d+1)t_0e^{2Lt_0^2}a_m.
$$
So $\p(\bar T=T_7)\le p_7$ for a suitably large $C$ by Chebyshev's inequality.
Set
\begin{align*}
\bar f(t)= & \sup_{k\le m-1}\sup_{s\le t}|Y_s^k|/\log(Na_k),\\
A^*_t= & \sup_{k\le m-1}\sup_{s\le t}|A^k_s|/\log(Na_k)
+\int_0^t|Y_s^m|ds/\log(Na_{m-1}),\\
D^*_t= & \sup_{k\le m-1}\sup_{s\le t}|D_s^k|/\log(Na_k).
\end{align*}
Since $Na_{m-1}>A\ge e$, we have $Na_k/\log(Na_k)\le Na_{k-1}/\log(Na_{k-1})$
for all $k\le m-1$. So we can use an argument from the proof of 
Proposition~\ref{QQ} to obtain $\bar f(t)\le e^{Lt}(A^*_t+D^*_t)$
for all $t\le t_0$.

The function $f(s)=\log s/(N^{1-d}\lambda s^d)$ is decreasing when $s>e$ and $e<A<Na_{m-1}\le
(AN^{d-1}/\lambda)^{1/d}$, so
$$
\log(Na_{m-1})/(Na_m)=f(Na_{m-1})\ge f((AN^{d-1}/\lambda)^{1/d})\ge\log N/(2A).
$$
Similarly, $\log(Na_{m-1})/\sqrt{Na_m}\ge\log N/(2\sqrt{A})$. Hence
$$
r(Na_m+\sqrt{Na_m})/\log(Na_{m-1})\le 4rA/\log N.
$$
We estimate as in the proof of Proposition~\ref{QQ} to obtain for $t<\bar T$
and $k\le m-1$
\begin{align*}
|b^k(X_t) & -b^k(x_t)-\nabla b^k(x_t)(X_t-x_t)|/\log(Na_k)\\
& \le H\lambda R^2(a_{k-1}^{d-1} +a_k^{d-1})/N\log(Na_k)\le 2HR^2/(N\log N)
\end{align*}
and 
$$
|Y_t^m|/\log(Na_{m-1}) \le(ra_m+N^{-1/2}r\sqrt{a_m})/\log(Na_{m-1})
\le 4rA/(N\log N).
$$
So $A^*_{\bar T}\le 2t_0(HR^2+2rA)/(N\log N) \le \half e^{-Lt_0}\bar R/N$,
provided $C$ is chosen suitably large.
For $k\le m-1$, set $\alpha_k=\frac12 e^{-Lt_0}\bar R\log(Na_k)/N$ and consider
the event $\O_k=\{\sup_{t\le\bar T}|D_t^k|>\alpha_k\}$.
On $\O_0=\O\setminus(\O_1\cup\dots\cup\O_{m-1})$ we have
$D^*_{\bar T}\le \half e^{-Lt_0}\bar R/N$, so $\bar f(\bar T)\le\bar R/N$
and $\bar T<T_8$.
Hence
$$
\p(\bar T=T_8)\le\sum_{k=1}^{m-1}\p(\O_k)
$$
and it will suffice to estimate $\p(\O_k)$ for each  $k\le m-1$.
Fix  $k\le m-1$. We can write $D_t^k=D_+(t)-D_-(t)$ where
$D_\pm(t)=D_\pm^1(t)+D_\pm^2(t)+D_\pm^3(t)$ and
\begin{align*}
D_\pm^1(t)=&N^{-1}(Z_\pm^k-W_\pm^k)(N\int_0^t\lambda_\pm^k(X_s)ds)\\
D_\pm^2(t)=&N^{-1}\{W_\pm^k(N\int_0^t\lambda_\pm^k(X_s)ds)
-W_\pm^k(N\int_0^t\lambda_\pm^k(\bar X_s)ds)\}\\
D_\pm^3(t)=&N^{-1/2}\int_0^t\{\s_\pm^k(\bar X_s)-\s_\pm^k(x_s)\}dB_\pm^k(s).
\end{align*}
Hence $\p(\O_k)\le q_+^1+q_-^1+q_+^2+q_-^2+q_+^3+q_-^3$ where,
for $j=1,2,3$, $q_\pm^j=\p(\O_\pm^j)$ and 
$\O_\pm^j=\{\sup_{t\le\bar T}|D_\pm^j(t)|>\alpha_k/6\}$.

For $t<\bar T$, we have $\lambda_\pm^k(X_t)\le\s^da_k$ so, taking
$t=N\s^da_kt_0$ and $x=N\a_k/6-c\log t$ in (\ref{KM}), we obtain
\begin{align*}
q_\pm^1&\le\p(\sup_{t\le N\s^da_kt_0}|Z(t)-W(t)|>N\a_k/6)\\
&\le c\s^dt_0Na_ke^{-N\a_k/(6c)}
=c\s^dt_0Na_ke^{-\bar R\log(Na_k)/(12ce^{Lt_0})}\\
&=c\s^dt_0(Na_k)^{1-\bar R/(12ce^{Lt_0})}
\le c\s^dt_0A^{1-\bar R/(12ce^{Lt_0})}\le C/(4A),
\end{align*}
for a suitable choice of $C$.

We turn to estimate $q^2_\pm$.
This will rely on the following continuity estimate
for Brownian motion:
for $\tau,h,\d>0$, setting $n=\lfloor\tau/h\rfloor$,
\begin{align*}
&\p(\sup_{s,t\le\tau,|s-t|\le h}|W(t)-W(s)|>\d)\\
&\quad\quad\le\p(\sup_{k\in\{0,1,\dots,n-1\},t\le2h}|W(kh+t)-W(kh)|>\d/2)\\
&\quad\quad\le2n\p(\sup_{t\le2h}W(t)>\d/2)\le(2\tau/h)e^{-\d^2/(16h)}.
\end{align*}
For $t<\bar T$ we have
\begin{align*}
\bar X^k_t\le&x^k_t+|X_t^k-x_t^k|+|\bar X_t^k-X_t^k|
\le\rho a_k+R\sqrt{a_k/N}+\bar R\log(Na_k)/N\\
\le&(\rho+R/\sqrt{A}+\bar R\log A/A)a_k\le\s a_k,
\end{align*}
so
$$
N\int_0^t\lambda_\pm^k(\bar X_s)ds\le\s^dt_0Na_k.
$$
Also, using the inequalities (\ref{B+}) and (\ref{B-}), for $t<\bar T$,
\begin{align*}
|N\int_0^t&(\lambda_+^k(X_s)-\lambda_+^k(\bar X_s))ds|\\
&\le Nd\s^{d-1}t_0\bar Ra_k\{\log(Na_{k-1})/(Na_{k-1})
+\log(Na_k)/(Na_k)\}\\
&\le2d\s^{d-1}t_0\bar R\log(Na_k)
\end{align*}
and
\begin{align*}
|N\int_0^t&(\lambda_-^k(X_s)-\lambda_-^k(\bar X_s))ds|\\
&\le\bar Rt_0\{\log(Na_k)+\log(Na_{k+1})1_{k\le m-2}\}
+rt_0(Na_m+\sqrt{Na_m})\d_{k,m-1}\\
&\le2t_0\bar R\log(Na_k),
\end{align*}
provided that $C$ is sufficiently large.
We take $\tau=\s^dt_0Na_k$, $h=2d\s^{d-1}t_0\bar R\log(Na_k)$
and 
$$
\d=N\a_k/6=\bar R\log(Na_k)/(12e^{Lt_0})
$$ 
to obtain
$$
q^2_\pm\le\frac{\s Na_k}{d\bar R\log(Na_k)}
\exp\{-\frac{\bar R\log(Na_k)}{4608d\s^{d-1}t_0e^{2Lt_0}}\}\le C/(4A)
$$
for a suitable choice of $C$.

It remains to estimate $q^3_\pm$. We shall show below that there exists a 
constant $C_0$ such that, for $t\le t_0$ and all $k\in\N$,
$$
\E(|\gamma^k_t|^2)\le C_0(x_t^{k-1}-x_t^k)\wedge(x_t^k-x_t^{k+1}).
$$
Then, for $t<\bar T$, 
\begin{align*}
&\sqrt N|\s_-^k(\bar X_t)-\s_-^k(x_t)|
=\sqrt N\left |\sqrt{((\bar X_t^k)^+-(\bar X_t^{k+1})^+)^+}-\sqrt{x_t^k-x_t^{k+1}}\right |\\
&\quad=\frac{\sqrt N|((\bar X_t^k)^+-(\bar X_t^{k+1})^+)^+-(x_t^k-x_t^{k+1})|}
                  {\sqrt{((\bar X_t^k)^+-(\bar X_t^{k+1})^+)^+}+\sqrt{x_t^k-x_t^{k+1}}}
\le(|\gamma_t^k|+|\gamma_t^{k+1}|)/\sqrt{x_t^k-x_t^{k+1}}
\end{align*}
so by Doob's $L^2$-inequality
$$
\E(\sup_{t\le\bar T}|D_-^3(t)|^2)
\le4\E\int_0^{t_0}N^{-1}|\s_-^k(\bar X_s)-\s_-^k(x_s)|^2ds\le16t_0C_0/N^2.
$$
Hence
$$
q_-^3=\p(\sup_{t\le\bar T}|D_-^3(t)|>\a_k/6)
\le2304t_0e^{2Lt_0}C_0/(\bar R\log(Na_k))^2\le C(\bar R\log A)^{-2}/2
$$
for a suitable choice of $C$.

The argument for $q^3_+$ is similar. For $t<\bar T$,
\begin{align*}
&\sqrt N|\s_+^k(\bar X_t)-\s_+^k(x_t)|\\
&\quad\quad=\sqrt N\left |\sqrt{(\lambda((\bar X_t^{k-1})^+)^d-\lambda((\bar X_t^k)^+)^d)^+}
-\sqrt{\lambda(x_t^{k-1})^d-\lambda(x_t^k)^d}\right |\\
&\quad\quad\le\lambda d\{(x_t^{k-1}\vee\bar X_t^{k-1})^{d-1}|\gamma_t^{k-1}|
+(x_t^k\vee \bar X_t^k)^{d-1}|\gamma_t^k|\}/\sqrt{\lambda(x_t^{k-1})^d-\lambda(x_t^k)^d}\\
&\quad\quad\le\lambda d2^{d-1}((x_t^{k-1})^{d-1}+(N^{-1/2}|\gamma_t^{k-1}|)^{d-1}) |\gamma_t^{k-1}|/\sqrt{\lambda(x_t^{k-1}-x_t^k)(x_t^{k-1})^{d-1}}\\
&\quad\quad\quad
+\lambda d2^{d-1}((x_t^k)^{d-1}+(N^{-1/2}|\gamma_t^k|)^{d-1})|\gamma_t^k|
/\sqrt{\lambda(x_t^{k-1}-x_t^k)(x_t^{k-1})^{d-1}}.
\end{align*}
so, by Doob's $L^2$-inequality,
\begin{align*}
\E(\sup_{t\le\bar T}|D_+^3(t)|^2)
& \le4\E\int_0^{t_0}N^{-1}|\s_+^k(\bar X_s)-\s_+^k(x_s)|^2ds\\
& \le8d^22^{2(d-1)}t_0(1+N^{(d-1)/2}C(d))C_0/N^2
\end{align*}
where $C(d)=\E(W(1)^{2d})$.  Hence
\begin{align*}
q_+^3
& =\p(\sup_{t\le\bar T}|D_+^3(t)|>\a_k/6)\\
& \le2152d^22^{2(d-1)}t_0(1+N^{(d-1)/2}C(d))C_0e^{2Lt_0}/(\bar R\log(Na_k))^2\\
&  \le C(\bar R\log A)^{-2}/2
\end{align*}
for a suitable choice of $C$. On combining this with the bounds for $q^1_\pm$ and $q^2_\pm$ 
already found, we obtain the desired bound for $p_8$.
\end{proof}

\begin{proof}[Proof of Theorem \ref{thm.diffusion}]
Set $A(N)=r(N)=(\log N)^{1/2}$ and 
define $\bar m(N)=\inf\{k\in\N: Na_k\le A(N)\}$.
Set $R(N)=(\log N)^{1/4}(1\wedge t_0)/2$.
It is straightforward to check that, if $C$ is the constant appearing in
Proposition~\ref{RR} and if $\bar R=3C$, then all the constraints in 
Propositions \ref{PP} and \ref{RR} are satisfied eventually and moreover
that $p_i(N)\to0$ for $i=1,2,3,4,7,8$. 
Since $\bar m(N)\ge m(N)$, this proves the theorem.
\end{proof}

\section{Fluctuation variance estimates}
\label{App}

We have deferred from other sections the analysis of certain linear
equations associated with our processes. The basic questions of existence 
and uniqueness in suitable spaces may be resolved by standard methods, so
we review this only briefly. The more delicate result, Proposition~\ref{GT},
which is needed for the diffusion approximation, relies on the particular
structure of our model.

We recall the $\nats\times\nats$ matrix-valued equation
\[\frac\partial{\partial t}\Phi_{t,s}=\nabla b(x_t)\Phi_{t,s},
\quad \Phi_{s,s}=I\]
to be solved for $0\le s\le t\le t_0$.
Note that, for $x_0\in S(\rho,t_0)$ and $t\in[0,t_0]$, we have 
$\|\nabla b(x_t)\|\le L$, where $\|\dots\|$ is the operator norm 
corresponding to $\|x\|=\sup_k|x^k|/a_k$.
Hence it is standard that this equation has a unique continuous 
solution with
$\|\Phi_{t,s}\|\le e^{L(t-s)}$ for all $s,t$.

The other relevant equations may be considered as stochastic perturbations
of the preceding equation. In Theorem~\ref{thm.jump.approx} we used the 
equation (\ref{eq.jump.approx})
$$
\tilde{\gamma}_t^k=\sqrt{N}\tilde M^k_t + \int_0^t \nabla b^k(x_s)
\tilde{\gamma}_s ds,\quad t\le t_0
$$
and in Theorem~\ref{thm.diffusion} we used the 
equation~(\ref{eq.diffusion})
\begin{equation}\label{GG}
\gamma^k_t=\sqrt{N}\bar 
M^k_t+ \int_0^t \nabla b^k(x_s) \gamma_s ds,\quad t\le t_0.
\end{equation}
Here 
$$
\tilde M^k_t=\int_{\mathbb R^{\mathbb N}\times(0,t]}
y^k(\tilde{\mu} -\tilde{\nu})(dy,ds)
$$
and
$$
\sqrt{N}\bar M^k_t=\int_0^t \sigma^k_+(x_s) dB^k_+(s) -\int_0^t \sigma^k_-(x_s) dB^k_-(s).
$$
Note that
\begin{align*}
N\E(|\tilde M^k_t|^2)
& =N\int_0^t\int_{\R^\N}(y^k)^2K(x_s,dy)ds\\
& =\int_0^t(\lambda_+^k(x_s)+\lambda_-^k(x_s))ds = N\E(|M^k_t|^2)
\end{align*}
and $\lambda_+^k(x_t)+\lambda_-^k(x_t)\le(\rho^d+1)a_k$ for $t\le t_0$.
Hence a standard type of iteration argument, using Doob's $
L^2$-inequality,
shows that these equations have unique measurable solutions with, respectively,
$$
\E(\sup_{t\le t_0}|\tilde\gamma_t^k|^2)\le8(\rho^d+1)t_0e^{2Lt_0^2}a_k
$$
and 
$$
\E(\sup_{t\le t_0}|\gamma_t^k|^2)\le8(\rho^d+1)t_0e^{2Lt_0^2}a_k.
$$
The details for equation~(\ref{eq.jump.approx}) follow below; 
equation~(\ref{eq.diffusion}) may be treated in the same way.
Note that for any vector $y$ and for all $k \in \nats$
\[ \frac{|\nabla b^k(x_s) y^k|}{\sqrt{a_k}} \le  L \sum_{j=k-1,k,k+1}
\frac{|y_{j}|}{\sqrt{a_j}}.\]
Now let
\[ \tilde \gamma_t^{(0)} = \sqrt{N} \tilde M_t ,\]
and for $n \in \nats$
\[ \tilde \gamma_t^{(n)}=\sqrt{N} \tilde M_t +  \int_0^t \nabla b^k(x_s)
\tilde \gamma_s^{(n-1)} ds.\]
Then for each $k$
\begin{eqnarray*}
\frac{ |\tilde \gamma^{(n+1),k}_t - \tilde \gamma^{(n),k}_t|^2 }{a_k} & \le &  
L^2 \left ( \int_0^t \sum_{j=k-1,k,k+1} \frac{ |\tilde \gamma^{(n),j}_s - 
\tilde \gamma^{(n-1),j}_s| }{\sqrt{a_j}} ds \right )^2\\
& \le & 3L^2 \sum_{j=k-1,k,k+1}  \left ( \int_0^t \frac{ |\tilde 
\gamma^{(n),j}_s - \tilde \gamma^{(n-1),j}_s| }{\sqrt{a_j}} ds \right )^2,
\end{eqnarray*}
so using Cauchy-Schwartz
\[ \frac{ |\tilde \gamma^{(n+1),k}_t - \tilde \gamma^{(n),k}_t|^2 }{a_k} \le 
3 L^2 t \sum_{j=k-1,k,k+1} \int_0^t  \frac{ |\tilde \gamma^{(n),j}_s - 
\tilde \gamma^{(n-1),j}_s|^2 }{a_j} ds.\]
Then for all $0 \le s \le t$
\[ \sup_{s \le t} \frac{ |\tilde \gamma^{(n+1),k}_s - \tilde 
\gamma^{(n),k}_s|^2 }{a_k} \le 
3 L^2 t \sum_{j=k-1,k,k+1} \int_0^t  \sup_{u \le s} \frac{ |\tilde 
\gamma^{(n),j}_u - \tilde \gamma^{(n-1),j}_u|^2 }{a_j} ds.\]
Let $h^{(n)} (t) = \sup_k \E ( \sup_{s \le t}|\tilde \gamma^{(n+1),k}_s - 
\tilde \gamma^{(n),k}_s|^2/a_k )$; then for $t \le t_0$
\[h^{(n)} (t) \le 9 L^2 t_0 \int_0^t h^{(n-1)} (s)ds.\]
Hence 
\[h^{(n)}(t) \le  \frac{2^{n} n!(3Lt_0)^{2n}}{(2n)!} 
\tilde M(t_0),\]
where
\begin{eqnarray*}
\tilde M (t)& = & \sup_k \E (\sup_{s\le t}{|\tilde M^k_s|}^2/a_k )\\
& \le & 4 \sup_k \sup_{s \le t} \E ({|\tilde M^k_s|}^2/a_k) \le 
4 t(\rho^d +1 ).
\end{eqnarray*} 
We deduce that $\tilde \gamma^{(n)}_t$ converges to a process $\tilde \gamma_t$
uniformly on $[0,t_0]$, and that $\tilde \gamma_t$ 
satisfies~(\ref{eq.jump.approx}). The uniqueness part of the proof is similar.
Let 
$$
g(t)=\sup_k\E(\sup_{s\le t}|\tilde \gamma_s^k|^2)/a_k.
$$
Then it follows from the above estimates that for $t \le t_0$
\[ g(t) \le \tilde M(t_0) e^{9L^2 t_0^2}.\]
%
%The reader who wishes to rederive these bounds, at least formally, may
%consider
%$$
%g(t)=\sup_k\E(\sup_{s\le t}|\gamma_s^k|^2)/a_k.
%$$
%Estimation using Doob's $L^2$-inequality leads to
%$$
%g(t)\le8(\rho^d+1)t_0+2t_0\int_0^tLg(s)ds.
%$$
%Assuming $g(t_0)<\infty$, Gronwall's lemma gives the claimed bounds. 
%For a proof, one takes this argument through the iteration scheme.
It may be verified by substitution that the formula (\ref{GF}) gives an
explicit representation of the solution of equation (\ref{eq.jump.approx}).

\begin{proposition}\label{GT}
The solution $(\gamma_t^k:k\in\nats,t\le t_0)$ to~(\ref{eq.diffusion}) 
satisfies
$$
\sup_{k\in\nats}\sup_{t\le t_0}\E(|\gamma_t^k|^2)
/\min\{x_t^{k-1}-x_t^k,x_t^k-x_t^{k+1}\}<\infty.
$$
\end{proposition}
\begin{proof}
Note that $\lambda_+^k(x_t)\le d(x_t^{k-1}-x_t^k)$, so
\begin{align*}
\dot x_t^k-\dot x_t^{k+1}
& =\lambda_+^k(x_t)-\lambda_+^{k+1}(x_t)
   -(x_t^k-x_t^{k+1})+(x_t^{k+1}-x_t^{k+2})\\
& \ge\lambda_+^k(x_t)-(d+1)(x_t^k-x_t^{k+1})+(x_t^{k+1}-x_t^{k+2}).
\end{align*}
Hence
$$
\int_0^t\lambda_+^k(x_s)ds+\int_0^t(x_s^{k+1}-x_s^{k+2})ds
\le e^{(d+1)t}(x_t^k-x_t^{k+1})
$$
and
$$
\int_0^t(x_s^k-x_s^{k+1})ds\le te^{(d+1)t}(x_t^k-x_t^{k+1}).
$$
Also
$$
\int_0^t\lambda_+^k(x_s)ds\le d\int_0^t(x_s^{k-1}-x_s^k)ds
\le dte^{(d+1)t}(x_t^{k-1}-x_t^k).
$$
Fix $\varepsilon>0$ and set
$$
f(t)=\sup_{k\in\nats}\sup_{s\le t}\E(|\gamma_s^k|^2)/\d_s^k
$$
where $\d_t^k=\min\{x_t^{k-1}-x_t^k,x_t^k-x_t^{k+1}\}+\varepsilon$.
Then $f(t_0)<\infty$.
Note that
\begin{align*}
&\E(|\int_0^t\sigma^k_+(x_s) dB^k_+(s)-\int_0^t \sigma^k_-(x_s) dB^k_-(s)|^2)\\
&\quad=\int_0^t\lambda_+^k(x_s)ds+\int_0^t(x_s^k-x_s^{k+1})ds
\le (dt+1)e^{(d+1)t}\d_t^k.
\end{align*}
We have
$$
\nabla b^k(x)\gamma=\lambda d(x^{k-1})^{d-1}\gamma^{k-1}
-(\lambda d(x^k)^{d-1}+1)\gamma^k+\gamma^{k+1}.
$$
We will make use of the following estimates
\begin{align*}
&\E(|\int_0^t\lambda d(x_s^{k-1})^{d-1}\gamma_s^{k-1}ds|^2)\\
&\quad\le\int_0^t\lambda^2d^2(x_s^{k-1})^{2(d-1)}\d_s^{k-1}ds
  \int_0^t\E(|\gamma_s^{k-1}|^2)/\d_s^{k-1}ds\\
&\quad\le d^2\int_0^t\lambda_+^k(x_s)ds\int_0^tf(s)ds
  \le (d^2+d^3t)e^{(d+1)t}\d_t^k\int_0^tf(s)ds
\end{align*}
and
\begin{align*}
&\E(|\int_0^t(\lambda d(x_s^k)^{d-1}+1)\gamma_s^kds|^2)\\
&\quad\le\int_0^t(\lambda d(x_s^k)^{d-1}+1)^2\d_s^kds
  \int_0^t\E(|\gamma_s^k|^2)/\d_s^kds
\le (d+1)^2te^{(d+1)t}\d_t^k\int_0^tf(s)ds
\end{align*}
and
$$
\E(|\int_0^t\gamma_s^{k+1}ds|^2)
\le\int_0^t\d_s^{k+1}ds
  \int_0^t\E(|\gamma_s^{k+1}|^2)/\d_s^{k+1}ds
\le e^{(d+1)t}\d_t^k\int_0^tf(s)ds.
$$
Now, from equation~(\ref{GG}), for all $t\le t_0$,
$$
\E(|\gamma_t^k|^2)\le A\d_t^k+B\d_t^k\int_0^tf(s)ds
$$
where $A=4(dt+1)e^{(d+1)t}$ and $B=8d(d+1)(dt+1)e^{(d+1)t}$. 
So $f(t)\le Ae^{Bt}$ by Gronwall's lemma. This bound 
does not depend on $\varepsilon$, so the proposition follows by letting
$\varepsilon\to0$.
\end{proof}


\begin{thebibliography}{10}


\bibitem{DN02} {\sc Darling, R.W.R. and Norris, J.R.} (2002). 
Structure of large random hypergraphs. To appear in {\it Ann. Appl. Prob.}


\bibitem{D77} {\sc Deimling, K.} (1977). 
\textit{Ordinary Differential Equations in Banach 
Spaces.} Lecture Notes in Mathematics 596, Springer-Verlag, Berlin.

\bibitem{ELZ86} {\sc Eager, D.L., Lazokwska, E.D. and Zahorjan, J.}
(1986). Adaptive load sharing in homogeneous distributed systems. {\it 
IEEE Trans. Soft. Eng.} \textbf{12} 662--675.

\bibitem{EK86}
{\sc Ethier, S.N. and Kurtz, T.G.} (1986). \textit{Markov Processes: 
Characterization and Convergence.} Wiley, New York.

\bibitem{G99}
{\sc Graham, C.} (2000). Kinetic limits for large communication networks. In
{\it Modelling in Applied Sciences} (N. Bellomo and M. Pulvirenti, 
eds.) 317--370. Birkh\"auser, Basel.

\bibitem{G00}
{\sc Graham, C.} (2000). Chaoticity on path space for a queueing network with 
selection of the shortest queue among several. {\em J. Appl. Prob.} {\bf 37}
198--201.

\bibitem{G04}
{\sc Graham, C.} (2004). Functional central limit theorems for a large network 
in which customers join the shortest of several queues. Preprint.


\bibitem{KMT75}
{\sc Koml\'os, J., Major, P. and Tusn\'ady, G.} (1975). 
An Approximation of Partial Sums of Independent RV's and the sample DF I. 
{\em Z. Wahrsch. verw. Gebiete} \textbf{32} 111--131.

\bibitem{L03}
{\sc Luczak, M.J.} (2003). A quantitative law of large numbers via exponential 
martingales. In {\it Stochastic Inequalities and Applications}
(E. Gin\'e, C. Houdr\'e and D. Nualart, eds.)
Birkh\"auser Series Progress in Probability {\bf 56} 93--111. Birkh\"auser, 
Basel.

\bibitem{LM03}
{\sc Luczak, M.J. and McDiarmid, C.} (2003). 
On the power of two choices: balls and bins in continuous time. Preprint.

\bibitem{LM04}
{\sc Luczak, M.J. and McDiarmid, C.} (2004). On the maximum queue length in 
the supermarket model. Preprint.


\bibitem{MS99}
{\sc Martin, J.B. and Suhov, Y.M.} (1998). Fast Jackson networks. {\it 
Ann. Appl. Prob.} \textbf{9} 854--870.

\bibitem{M96}
{\sc Mitzenmacher, M.} (1996). The power of two choices in randomized load 
balancing. PhD thesis, Berkeley.


\bibitem{RY91}
{\sc Revuz, D. and Yor, M.} (1991). \textit{Continuous Martingales and 
Brownian Motion.} Springer-Verlag, Berlin.

\bibitem{VDK96}
{\sc Vvedenskaya, N.D., Dobrushin, R.L. and Karpelevich, F.I.} (1996). 
Queueing system with selection of the shortest of two queues: an asymptotic 
approach. {\em Prob. Inf. Transm.} \textbf{32} 15--27.

\end{thebibliography}
\end{document}